# Examining Gender Differences in the Impact of Instructional Clarity on Year 9 Students' Mathematics Interest in New Zealand: A Multi-Group Comparison of 2019 TIMSS Data

Huayu Gao
*University of Auckland*
hgao462@aucklanduni.ac.nz

Tanya Evans
*University of Auckland*
t.evans@auckland.ac.nz

Gavin T.L. Brown
*University of Auckland*
gt.brown@auckland.ac.nz

Grounded in Social Cognitive Career Theory (SCCT), this study explores how teacher instruction clarity impacts Year 9 students' mathematics interest through the mediating roles of confidence and value, with a focus on gender differences. Structural equation modelling applied to the 2019 New Zealand TIMSS dataset (2,594 girls and 2,651 boys) supports the SCCT framework. While confidence and value have similar mediating effects for boys, confidence plays a more pivotal role for girls. Gender differences are primarily observed in the value-mediated pathway, underscoring persistent challenges to achieving educational equity through teachers' instruction in New Zealand.

Gender disparities in New Zealand's mathematics-intensive fields remain a persistent issue. While the proportion of female graduates in NZ tertiary education has reached an impressive 61.7%, surpassing that of their male counterparts, women's representation in STEM (Science, Technology, Engineering, and Mathematics) fields remains disproportionately low at only 39.9% (Organisation for Economic Co-operation and Development [OECD], 2023). This gender segregation trend further extends to the occupation domain, as exemplified by the University of Auckland (NZ), where female faculty members in the Faculty of Science and Engineering are significantly underrepresented compared to their male colleagues, particularly in engineering disciplines (Brower & James, 2020). Without effective intervention, the current gender imbalance may persist until 2070 (Brower & James, 2020), hindering the advancement of gender equality and posing severe challenges to the cultivation of a diverse talent pool.

Social Cognitive Career Theory (SCCT) offers a comprehensive theoretical framework for explaining gender differences in career choices. The theory posits that interest serves as a crucial source for individuals to set career goals; thus, investigating the formation mechanism of interest is paramount in unveiling gender disparities in the career decision-making process. According to the interest model of SCCT, self-efficacy and outcome expectations are the core social cognitive variables that influence interest formation. Self-efficacy pertains to an individual's belief in their capability to accomplish specific tasks, whereas outcome expectations refer to an individual's anticipation of the potential outcomes of particular behaviours and their assessment of the value associated with these outcomes (Lent et al., 1994). Bandura (1977) emphasised that self-efficacy typically wields a more potent influence on behavioural decisions. Furthermore, given that individuals tend to anticipate ideal outcomes in activities they perceive themselves to be proficient in, outcome expectations are partially determined by self-efficacy. Contextual and experiential factors, as antecedents of these social cognitive beliefs, exert direct or indirect effects on interest formation (Lent et al., 1994). Specifically, these extra-personal factors not only provide information and resources related to career development but also shape individuals' beliefs, thereby influencing the trajectory of interest development. In the process of interest formation, gender is considered a key moderating variable. SCCT posits that gender is not merely a biological attribute but also a characteristic formed through social construction. The socialisation of gender roles renders individuals' exposure to career-related experiences and information selective, resulting in divergent paths of influence for boys and girls in the interest formation process.

Many empirical studies support the explanations of the SCCT's interest model. Research indicates that girls' choices of STEM pathways are significantly influenced by their interest in these fields (Blotnicky et al., 2018). Moreover, the development of individuals' interest in STEM careers is directly affected by their self-efficacy and values regarding STEM occupations and indirectly influenced by environmental factors such as family, peers, and school (Turner et al., 2019). For instance, Wang et al. (2023) validated the applicability of the interest model in SCCT within the Chinese cultural context, demonstrating that self-efficacy and individuals' perceptions of career prospects, required skills, and self-development significantly enhance their interest in STEM fields. They also found that formal or informal learning experiences and social support indirectly influence the development of interest by enhancing self-efficacy and perceptions of STEM careers. Additionally, their study highlighted the role of sociocultural factors in the gender differences observed in interest development, suggesting that self-efficacy may play a more important role in the formation of girls' interest in STEM, particularly in cultural backgrounds where gender role stereotypes are more prominent.

School factors, as framed by SCCT, have increasingly become a focal point in educational research. Given that adolescents spend a substantial amount of time in school, teachers as critical agents between students and schools, significantly influence students' learning experience and motivation. Research demonstrates that teachers have the capacity to foster positive perceptions of STEM fields among girls (Smith & Evans, 2024), and that the quality of teachers' instruction is closely associated with students' learning motivation and academic performance (Ekmekci & Serrano, 2022). Moreover, effective pedagogical strategies provided by educators create conducive learning environments that enhance girls' mathematical interests (Kang & Keinonen, 2017).

Instructional clarity has been recognised as a crucial aspect of effective teaching, particularly in mathematics, where clear instruction is linked to variations in students' learning outcomes. For example, a comparative analysis of TIMSS 2019 data from Hong Kong and England was conducted by Chen and Lu (2022), revealing that instructional clarity not only improves students' enjoyment and achievement in mathematics but also reduces boredom in both jurisdictions. Likewise, Zhang et al. (2025) utilised 2019 Australia TIMSS data and discovered a significant positive correlation between instructional clarity and students' mathematics achievement, and further demonstrated that teacher instructional clarity could positively predict learning motivation. Moreover, Titsworth et al. (2015) conducted two meta-analyses, which highlighted the overall moderate positive impact of teachers' instructional clarity on students' affective learning. Nevertheless, the high heterogeneity observed in Titsworth et al.'s (2015) results suggested that contextual or individual factors may moderate the impact of instructional clarity. Notably, gender is widely acknowledged as a key moderating variable in TIMSS studies, with existing research consistently demonstrating the systematic influence of gender differences on students' learning attitudes and motivation, particularly with boys exhibiting higher confidence and intrinsic values in mathematics (Michaelides et al., 2019; Watt et al., 2012).

Building on these insights, the present study adopts the SCCT framework to examine the role of instructional clarity in fostering students' interest in mathematics. It further investigates whether there are gender differences in the pathways through which teachers' instructional clarity influences students' mathematics interests. By elucidating the critical role of instructional clarity in shaping learning motivation and promoting gender equity, this study aims to provide empirical evidence for optimising mathematics teaching practices.

Additionally, it holds the potential to uncover the underlying mechanisms that drive students' academic choices.

## Research Questions

Guided by the SCCT, this study constructs a hypothetical model, as shown in Figure 1 and focuses on the following three core research questions:

- What are the specific mechanisms through which instruction clarity influences students' mathematics interests via confidence and values?
- Among the three mediating paths shown in Figure 1 (Instruction Clarity → Confidence → Interest; Instruction Clarity → Value → Interest; Instruction Clarity → Confidence → Value → Interest), which path has the most significant impact on boys' mathematics interest? Which path has the most significant impact on girls' mathematics interests?
- Are there significant gender differences in the three mediating paths shown in Figure 1? If such differences exist, in which mediating path are they most pronounced?

**Figure 1**

*The Hypothetical Model Illustrating the Relation between Instruction Clarity and Mathematics Interest*

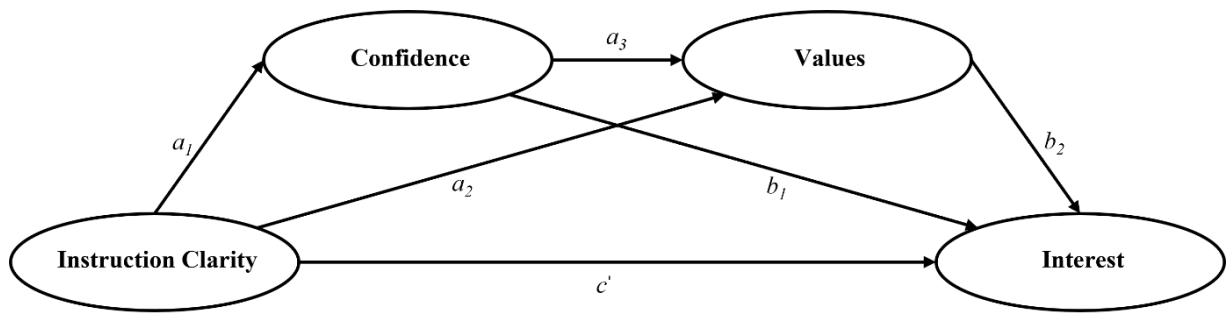

## Methodology

### Participants

The study comprised 2,594 Year 9 girls and 2,651 Year 9 boys in the first year of New Zealand secondary education from the publicly available 2019 TIMSS dataset. TIMSS 2019 collected contextual information through the administration of questionnaires to students, their teachers, and principals.

### Instruments

The 2019 TIMSS background questionnaires collected data on students' attitudes and motivation toward mathematics study, employing a theoretical framework grounded in psychological constructs, including latent variables such as confidence, interest, value, and instructional clarity which students perceive in mathematics classes. Specifically, instructional clarity was measured through five four-point Likert-type items in the 2019 TIMSS, assessing the clarity of instruction provided by teachers during mathematics lessons from students' perspectives, which were aggregated into a scale termed "Instructional Clarity" in this study. Students' confidence in mathematics was measured via six four-point Likert-type items, constituting a scale variable termed "Mathematics Confidence". Six questionnaire items captured students' mathematical values, forming a scale variable termed "Mathematics Value". Additionally, interest in mathematics was assessed through seven four-point Likert-type items evaluating students' enjoyment in learning mathematics, which were combined to create the scale variable termed "Mathematics Interest". To facilitate data processing, reverse scoring was applied to some items, ensuring that higher scores reflected higher levels of the measured items.

These four selected scales from the 2019 TIMSS can be interpreted within the SCCT framework, and their reliability is described in detail below:

The *Instruction Clarity (IC)* scale consists of five items: a. I know what my teacher expects me to do, b. My teacher is easy to understand, c. My teacher has clear answers to my questions, d. My teacher is good at explaining mathematics, e. My teacher explains a topic again when we don't understand. The internal consistency reliability of the scale was assessed using Cronbach's α coefficient, yielding a result of 0.895, indicating high reliability.

The *Mathematics Confidence (MC)* scale comprises six items: a. Mathematics is more difficult for me than for many of my classmates, b. Mathematics is not one of my strengths, c. Mathematics makes me nervous, d. My teacher tells me I am good at mathematics, e. Mathematics is harder for me than any other subject, f. Mathematics makes me confused. The Cronbach's α for this scale was 0.828, indicating good reliability.

The *Mathematics Value (MV)* scale includes six items: a. I need mathematics to learn other school subjects, b. I need to do well in mathematics to get into the university of my choice, c. I would like a job that involves using mathematics, d. It is important to learn about mathematics to get ahead in the world, e. Learning mathematics will give me more job opportunities when I am an adult, f. My parents think that it is important that I do well in mathematics. The Cronbach's α for this scale was 0.848, indicating good reliability.

The *Mathematics Interest (MI)* scale consists of seven items: a. I enjoy learning mathematics, b. I wish I did not have to study mathematics, c. I learn many interesting things in mathematics, d. I like mathematics, e. I like to solve mathematics problems, f. I look forward to mathematics class, g. Mathematics is one of my favourite subjects. The Cronbach's α for this scale was 0.930, indicating good reliability.

**Method**

To investigate the relations among the four latent variables, "Instruction Clarity" (IC), "Mathematics Confidence" (MC), "Mathematics Values" (MV), and "Mathematics Interest" (MI), this study employed a multiple mediator model within the framework of structural equation model using SPSS 30.0 and AMOS 30.0 software. This method offers several notable advantages (Fang et al., 2014). First, it enables the examination of the specific mediation effect of each mediator variable while controlling for other mediators, such as $a_1b_1$, $a_2b_2$, and $a_1a_3b_2$. Second, it provides contrast mediation effects, such as $a_1a_3b_2 - a_2b_2$, $a_1b_1 - a_2b_2$, and $a_1a_3b_2 - a_1b_1$, allowing researchers to determine which mediator variable has a stronger influence.

## Results

The study assessed model fit across genders using separate samples, yielding acceptable fit indices (For boys: $\chi^2$ = 2700.056, RMSEA = 0.060, CFI = 0.940; For girls: $\chi^2$ = 2796.916, RMSEA = 0.060, CFI = 0.930) that met conventional standards (West et al., 2012). A smaller $\chi^2$ value suggests a better model fit. The RMSEA measures the approximate error of the model fit, with values below 0.08 indicating a good fit. The CFI compares the model to a baseline (worst-fit) model, with values closer to 1 indicating a better fit.

Subsequently, a cross-gender measurement invariance analysis was conducted to examine the consistency of the questionnaire structure across different gender groups. Given this study did not involve mean or variance comparison, the following sequence of tests was performed: (1) *Configural Invariance*, ensuring that different groups have the same factor structure, including consistent factor numbers and factor-item correspondence; (2) *Metric Invariance* or *Weak Invariance*, testing whether factor loadings are invariant across groups. The configural invariance test results for boys and girls showed acceptable levels of RMSEA (0.044) and CFI (0.934), supporting configural invariance. The metric invariance test results supported metric invariance based on traditional criteria, with ΔRMSEA = 0.001 and ΔCFI < 0.001, both less than 0.01, and $\Delta\chi^2$ = 61.400, $\Delta df$ = 20, $p$ = 0.070, which was not significant.

After confirming the good fit of the SEM, the significance of the mediation effect was tested using the bias-corrected bootstrap method. If the confidence interval of a mediation effect does not include 0, the effect is considered significant (Fang et al., 2014). It is important to note that a significant mediation effect only indicates that it is not equal to 0, and further analysis of the mediation effect size requires comparing the mediation effects of different paths.

The analysis results showed that "IC" had a significant positive impact on students' "MI" (boys: $c'$ = 0.260, $p$ < 0.001; girls: $c'$ = 0.270, $p$ < 0.001). Additionally, IC influenced MI through three mediation paths formed by "MC" and "MV" (see Figure 1 & Table 1). For boys, the following three mediation paths were significant: IC → MC → MI: $a_1b_1$ = 0.152; z = 8.444; 95% CI [0.117, 0.188], IC → MV → MI: $a_2b_2$ = 0.166; z = 11.067; 95% CI [0.139, 0.202], IC → MC → MV → MI: $a_1a_3b_2$ = 0.017; z = 5.667; 95% CI [0.011, 0.025]. For girls, the three paths were also significant: IC → MC → MI: $a_1b_1$ = 0.191; z = 11.235; 95% CI [0.161, 0.229], IC → MV → MI: $a_2b_2$ = 0.108; z = 8.308; 95% CI [0.085, 0.136], IC → MC → MV → MI: $a_1a_3b_2$ = 0.022; z = 5.500; 95% CI [0.016, 0.031]. When comparing the mediation effects, it was found that the difference between the mediation effects of MC and MV was not significant (z = 0.560, 95% CI [−0.033, 0.066]) among boys. However, in the girls' group, the mediation effect of MC was significantly greater than that of MV (z = −3.609, 95% CI [−0.128, −0.038]). Furthermore, in the gender difference analysis, only the mediation effect of MV showed a significant gender difference, with the mediation effect of MV being significantly greater among boys compared to girls (z = −2.900, 95% CI [−0.101, −0.019]).

**Table 1**

*SEM Results of Multiple Mediation Effects of Instruction Clarity on Mathematics Interest*

| Gender | Mediation paths | Point estimate | SE | z | Bias-corrected percentile | |
|---|---|---|---|---|---|---|
| | | | | | Lower | Upper |
| Boys | IC→MC→MI | 0.152 | 0.018 | 8.444 | 0.117 | 0.188 |
| | IC→MV→MI | 0.166 | 0.015 | 11.067 | 0.139 | 0.202 |
| | IC→MC→MV→MI | 0.017 | 0.003 | 5.667 | 0.011 | 0.025 |
| Girls | IC→MC→MI | 0.191 | 0.017 | 11.235 | 0.161 | 0.229 |
| | IC→MV→MI | 0.108 | 0.013 | 8.308 | 0.085 | 0.136 |
| | IC→MC→MV→MI | 0.022 | 0.004 | 5.500 | 0.016 | 0.031 |
| | Contrast Mediation | | | | | |
| Boys | MC − MCMV | 0.135 | 0.016 | 8.438 | 0.104 | 0.169 |
| | MV − MC | 0.014 | 0.025 | 0.56 | −0.033 | 0.066 |
| | MV − MCMV | 0.149 | 0.016 | 9.313 | 0.122 | 0.184 |
| Girls | MC − MCMV | 0.169 | 0.016 | 10.563 | 0.142 | 0.205 |
| | MV − MCMV | 0.086 | 0.014 | 6.143 | 0.06 | 0.115 |
| | MV − MC | −0.083 | 0.023 | −3.609 | −0.128 | −0.038 |
| Girls − Boys | MC | 0.040 | 0.025 | 1.600 | −0.011 | 0.086 |
| | MV | −0.058 | 0.020 | −2.900 | −0.101 | −0.019 |
| | MCMV | 0.005 | 0.005 | 1.000 | −0.005 | 0.015 |

## Discussion

### The Role of Instructional Clarity in Enhancing Students' Mathematics Interest

The present study revealed that in the New Zealand context, teacher instructional clarity not only directly enhances students' interest in mathematics but also indirectly fosters their interest by increasing their confidence and perceived value of mathematics. Notably, the improvement in the perceived value of mathematics is, to some extent, dependent on the enhancement of confidence. This finding supports the theoretical assumptions of SCCT regarding the mechanisms of interest development and is consistent with the conclusions of the majority of existing empirical studies (Wang et al., 2023). However, Cheung and Sonkqayi (2025) suggested that in Asian regions, instructional clarity exerts negative effects on students' science performance. The conflicting conclusions about the impact of teacher instructional clarity on the students' learning outcomes may be attributable to differing teaching styles shaped by distinct sociocultural contexts. In Asian regions, the teacher-centred instructional approach is prevalent, enabling educators to guide students regularly in solving academic problems. In contrast, countries such as New Zealand and Australia advocate for a student-centred approach, which typically involves less teacher intervention. Moreover, Niu et al. (2025) demonstrated that in countries with student-centred learning environments, teacher guidance significantly enhances students' academic performance compared to those in teacher-centred education systems. Therefore, within the educational setting of New Zealand, teachers' clear instruction and support may be more conducive to enhancing students' interest and achievement. It is essential that mathematics teachers in New Zealand should actively engage in relevant training programs and seek professional development support to provide students with more structured and clear guidance, which could effectively promote their interest in mathematics.

### The Gender-Specific Mediating Patterns in the Development of Mathematics Interest

The present study further uncovered gender-differentiated mediation mechanisms, whereby confidence and values contribute in distinct ways to the development of interest. Specifically, confidence and values exhibit similar mediating effects in fostering interest among boys. However, for girls, the mediating effect of confidence is significantly stronger than that of values. These findings in our study are in accord with the majority of previous research, indicating that confidence plays a more important role in girls' mathematical interests, especially in sociocultural contexts with strong gender role stereotypes (Wang et al., 2023). Traditional gender role stereotypes often frame mathematics as a "male-dominated domain", which may unconsciously lead teachers to have lower expectations for girls in mathematics (Smith & Evans, 2024). Research also indicates that in coeducational learning environments, girls interact with teachers less frequently than boys (Smith & Evans, 2024). Given girls' heightened sensitivity to social interaction cues in the classroom (Reinholz et al., 2022), this disparity in interaction frequency makes girls more likely to perceive and internalise teachers' lower expectations in mathematics classes. However, confidence in mathematics may help girls develop intrinsic beliefs in their abilities, which serves as a buffer against the negative effects of gender role stereotypes or lower expectations from educators. This confidence strengthens their self-identity, enabling them to sustain their enthusiasm and commitment to mathematics even when faced with societal biases. Considering the importance of confidence in nurturing girls' interest in mathematics, teachers should be mindful of girls' psychological needs in the classroom. By creating a positive and inclusive learning environment, providing constructive feedback and encouragement, teachers can effectively cultivate girls' confidence and further stimulate their interest in mathematics.

## The Gender-Differentiated Mediating Pathways in the Development of Mathematics Interest

While boys and girls demonstrate different mediating mechanisms in shaping their interest in mathematics, our study revealed that the most salient gender-specific distinction lies in their relative dependence on values rather than confidence. Specifically, the mediating role of confidence between instructional clarity and interest remains consistent across genders. However, the mediating role of values is markedly more pronounced among boys. The similar weights of confidence as a mediating factor between boys and girls suggest the relative stability of confidence in fostering mathematics interest, aligning with Bandura's (1977) assertion that robust confidence is a uniformly accurate predictor of performance under varying challenges. Conversely, the observed gender differences in the mediating effects of values on students' interest in mathematics imply that teachers' instruction may inadvertently perpetuate educational inequities in New Zealand. This finding is consistent with evidence from studies on high school students in the United States, where a disconnect between girls' values and interests in mathematics has been observed, suggesting that girls may rely less on values than boys in the development of mathematical interest (Wang et al., 2015). A plausible explanation for this gender-specific disparity in the development of interest may reside in teachers' implicit gender role beliefs, which may inadvertently foster more supportive mathematical environments for boys, thereby reinforcing the mediating role of value in shaping boys' mathematics interests. These findings in the present study highlight the potential limitations of teachers' instructional strategies in mitigating gender disparities in the formation of mathematics interests. As such, it is imperative that educators recognise and address gender differences through targeted interventions during the instructional process in New Zealand.

## Conclusion

This study employed a multiple mediation model to investigate the mediating roles of confidence and value in the relationship between teacher instructional clarity and students' mathematical interest among boys and girls in their first year of secondary school in New Zealand. The findings revealed that clear instruction from teachers not only directly enhances students' interest in mathematics but also indirectly fosters their interest through the mediating effects of confidence and value. Notably, confidence emerges as a critical mediator for girls, while the primary gender difference lies in the mediating effect of value. Specifically, the value boys attach to mathematics appears to play a more prominent role in fostering their interest than it does for girls. These results highlight that, while teacher instructional clarity can effectively promote students' interest in mathematics, gender role stereotypes embedded in teachers' expectations may be transmitted to students through classroom interactions, thereby influencing students' beliefs and undermining educational equity. However, it is imperative to recognise that the influence of gender roles may vary across different sociocultural contexts. The findings of this study are specific to the New Zealand context, and further empirical research is required to confirm their generalisability.

## Acknowledgement

Ethics approval Ref. UAHPEC24842 was granted by the University of Auckland.